% AMS-TeX,  amsppt style

\documentstyle{amsppt}
\loadbold
\NoBlackBoxes
%\magnification=1200

\hsize = 5 true in
\vsize = 7.5 true in

\hoffset .8 true in
\voffset .8 true in
\vskip 1.4 true in

\define\eqdef{\overset \text{def}\to =}
\define\LL{\Cal L}

\define\ZZ{\Bbb Z}
\define\FF{\Bbb F}

\define\PP{\Bbb P}

\define\lam{{\lambda}}
\define\Lam{{\Lambda}}
\define\Del{\Delta}
\define\del{\delta}
\define\om{{\varpi}}

\define\ep{{\epsilon}}
\define\ka{{\kappa}}
\define\tka{\widetilde{\kappa}}
\define\ii{\bold i}

%\define\tildii{  \tilde{\boldsymbol{\imath}} }
\define\tJ{\widetilde{J}}

\define\ttau{\tilde{\tau}}
\define\tii{  \widetilde{\bold i} }
\define\mm{\bold m}
\define\Zii{{Z_{\ii}}}

\define\OO{{\Cal O}}
\define\TT{{\Cal T}}
\define\CC{{\Cal C}}
\define\Lm{{\LL_{\mm}}}

\define\Gr{\roman{Gr}}

\define\tr{\roman{tr}}
\define\Char{\roman{char}^*\, }
\define\diag{\roman{diag}}

\define\T{\Cal T}

\define\nul{{\bold O}}
\define\st{\star}

\define\sq{\times}
\define\emp{\emptyset}
\define\miss{\circ}

\define\wmax{w_{\roman{max}}}
\define\kmax{\ka_{\roman{max}}}
\define\tkmax{\tka_{\roman{max}}}
\define\wJ{w_J}
\define\lcol{\overset k\to <}

\define\forsome{\exists}

\topmatter

\title 
Standard Monomial Theory for Bott-Samelson Varieties
of $GL(n)$
\endtitle

\author V. Lakshmibai and Peter Magyar \endauthor

\address{ Department of Mathematics, 
Northeastern University, Boston, MA 02115, USA}
\endaddress

\email lakshmibai\@neu.edu \quad
pmagyar\@lynx.neu.edu \endemail

\thanks Both authors partially supported by the 
National Science Foundation. \endthanks

\abstract
We construct an explicit basis for the 
coordinate ring of the Bott-Samelson variety $\Zii$
associated to $G = GL(n)$ and an arbitrary
sequence of simple reflections $\ii$.  Our basis is
parametrized by certain standard tableaux and generalizes
the Standard Monomial basis for Schubert varieties.
\endabstract

%\leftheadtext{ Standard Monomial Theory}
\rightheadtext{ Standard Monomial Theory}

\endkeywords

%\subjclass {17B20,17B70}
%\endsubjclass

\endtopmatter

\document

In this paper, we prove 
the results announced in [LkMg] 
for the case of Type {\bf A$_{n-1}$} (the groups
$GL(n)$ and $SL(n)$).
That is, we construct an explicit
basis for certain ``generalized
Demazure modules'', natural finite-dimensional
representations of the group $B$ of upper triangular matrices.  
These modules can be constructed in an elementary way
as flagged Schur modules [Mg1,Mg2,RS1,RS2].
They include as special cases almost all
natural examples of $B$-modules, and their characters
include most of the known generalizations of Schur polynomials.
We view these representations geometrically via
Borel-Weil theory as the space of global sections 
of a line bundle over a  Bott-Samelson variety.  
Thus, our theory also describes the coordinate
ring of this variety.

Notations: $G = GL(n,\FF)$ or $SL(n, \FF)$,  
where $\FF$ is an algebraically closed field of 
arbitrary characteristic or $\FF = \ZZ$\,;  
\  $B$ is the Borel subgroup consisting of upper triangular
matrices;\ \, $T$ is the maximal torus consisting of diagonal matrices;\ 
$W$ is the symmetric group $S_n$
generated by the adjacent transpositions 
(simple reflections)
$s_i=(i,i+1)$;\ \
$P_i \supset B$ is
the minimal parabolic subgroup of $G$ 
associated to $s_i$, namely
$P_i = 
\{\, (x_{ij}) \in G \mid x_{ij}=0 \text{ if }
   i>j \text{ and } (i,j)\! \neq\! (i\!+\!1,i)\, \}.$ 

For any word $\ii = (i_1,\ldots,i_l)$, with letters
$1\leq i_j \leq n-1$,
the {\it Bott-Samelson variety} is the quotient space
$$
\Zii = P_{i_1} \times P_{i_2} \times \cdots \times P_{i_l} 
\ / \ B^l \ ,
$$
where $B^l$ acts by 
$$
(p_1, p_2, \ldots, p_l)\cdot (b_1,b_2,\ldots,b_l) = 
(p_1 b_1, b_1^{-1} p_2 b_2,\ldots, b_{l-1}^{-1} p_l b_l).
$$
It was originally used  [BS, D1] to desingularize
the Schubert variety $X_w = \overline{B\cdot w B}
\subset G/B$, where
$w = s_{i_1} \cdots s_{i_l}$.
The desigularization is given by the multiplication map 
$\Zii \rightarrow X_w \subseteq G/B$,\
$(p_1,\ldots,p_l) \mapsto p_1\!\cdots\! p_l\!\cdot\! B$,
and $\Zii$ has the structure of an iterated fiber 
bundle with fiber $\PP^1$ in each iteration, 
so we may loosely think of $\Zii$ as a ``factoring'' of the
Schubert variety into a twisted product of projective lines.

Denote $\Gr(i) = \Gr(i,\FF^n)$ 
the Grassmannian of 
$i$-dimensional subspaces of linear $n$-space, and  
$$
\Gr(\ii) \eqdef \Gr(i_1) \times \cdots \times \Gr(i_l).
$$ 
We can realize $\Zii$ as a variety of configurations
of subspaces of $\FF^n$ 
(a kind of multiple Schubert variety)
via the embedding by successive multiplications [Mg2]:
$$
\matrix 
\mu :\,& \Zii & \rightarrow & \Gr(\ii) \\
&(p_1,\ldots,p_l) & \mapsto & 
(p_1 \FF^{i_1},\, p_1 p_2 \FF^{i_2}, \ldots,\,
p_1\!\! \cdots\! p_l \FF^{i_l})
\endmatrix
$$
where $0 \subset \FF^1 \subset \cdots \subset \FF^n$ 
is the standard flag fixed by $B$.
Although we will not need it here,
we note that $\mu(\Zii)\cong \Zii$ can be described explicitly in terms of 
incidence relations: that is, a configuration of subspaces
$(V_1,\ldots,V_l) \in \Gr(\ii)$ lies in $\mu(\Zii)$ exactly if
certain inclusions $V_i\subset V_j$ are satisfied, 
as specified by the combinatorics of wiring diagrams.
See [Mg2].

Now, each Grassmannian has a 
minimal-degree ample line bundle 
(the Plucker bundle) $\OO(1)$,
and for any sequence $\mm=(m_1,\ldots,m_l)$,
$\ m_j \in \ZZ_+$, there is an effective line bundle
on $\Gr(\ii)$ given by tensoring the $m_j$th
powers of the Plucker bundles on the factors of $\Gr(\ii)$:
$
\OO(\mm) = \OO(1)^{\otimes m_1} \otimes 
\cdots \otimes \OO(1)^{\otimes m_l}.
$
Denote its restriction to $\Zii$ by
$\Lm = \mu^* \OO(\mm)$.
We shall be concerned with the $B$-module
$$
H^0(\Zii,\Lm),
$$
which includes as special cases
dual Schur modules (Weyl modules) [FH,F], 
Demazure modules [D1,LkSh2,LkSb1], skew Schur modules [FH,F],
the Schubert modules of Kraskiewicz and Pragacz [KP],
and the generalized Schur modules of
percent-avoiding diagrams [RS2].
Thus, the characters of these modules include Schur,
key, skew Schur, and Schubert polynomials.  See [Mg2].
In particular, we obtain a new proof of the classical
Standard Monomial Theory for Demazure modules of type $A$.

An example will give the flavor of our results.
Let $G= GL(3)$, $\ \ii = (1,2,1)$, $\ \mm=(1,1,1)$.  
We may write
$$
\Gr(\ii) = \Gr(1) \times  \Gr(2) \times  \Gr(1)
$$
$$
\Zii \cong \{\ (V_1,V_2,V'_1) \in 
\Gr(\ii) \ \mid  \
\FF^2 \supset V_1 \subset V_2 \supset V'_1 \ \}.
$$
Then $H^0(\Zii,\Lm)$ is spanned by
all products of the form 
$$
\Del_{abcd} = \Del_a(x) \Del_{bc}(y) \Del_d(z),
$$
where $1\leq a,b,c,d\leq 3$ and $\Del_{a}$, $\Del_{bc}$, $\Del_d$ 
mean minors on the corresponding rows of the homogeneous
coordinates on $\Gr(\ii)$:
$$
(x,y,z)\, = 
\pmatrix x_1 \\ x_2 \\ x_3 \endpmatrix\!\! \times \!\!
\pmatrix y_{11} & y_{12} \\ y_{21} & y_{22} \\ y_{31} & y_{32}
\endpmatrix \!\!\times\!\!
\pmatrix z_1 \\ z_2 \\ z_3 \endpmatrix
\in \Gr(\ii) .
$$
For example, $\Del_{2132} = x_1(y_{11} y_{32}-y_{31}y_{12})z_2$.
The sequence of indices $\tau = abcd$ 
indexing a spanning vector $\Del_{abcd}$
is called a {\it tableau}.
Theorems 1 and 2 below allow us to select a basis 
of $H^0(\Zii,\Lm)$ from
the spanning set, corresponding to the set of {\it standard
tableaux}:
$$
\split
\tau \in \TT(\ii,\mm) = \{ 1 12 1, 2 12 1, 2 12 2 ,
1 13 1, 2 13 1, 2 23 1, 2 23 2,\\ 
1 12 2,
1 13 2, 2 13 2,
 1 13 3, 2 13 3, 2 23 3 \}  \quad 
\endsplit
$$
Since each $\Del_{abcd}$ is an eigenvector of the
diagonal matrices, this allows us to compute the
character of the $B$-module $H^0(\Zii,\Lm)$.

In the general case, we give two descriptions of the 
standard tableaux $\TT(\ii,\mm)$.  The first (\S 1.2)
is in the spirit of the monotone lifting property 
of classical Standard Monomial Theory [LkSd1, LkSd2] 
(which in turn generalizes 
Young's increasing-in-rows-and-columns definition).
The second (\S1.4) is in terms of the refined Demazure
character formula and 
crystal lowering operators of Lascoux and Schutzenberger
[LcSb1] and Littelmann [Lt1, Lt2, Lt3].  This
description is more suited to computations,
and it gives an efficient algorithm for generating 
the standard tableaux.
The above list of 13 tableaux, for example, can be 
computed by hand in less than a minute.

The paper is organized as follows.
In Section 1 we define the standard tableaux
and state the main theorems.
In Section 2 we prove the equivalence of our
two definitions of standard tableaux by an
elementary argument.  In Section 3 we show
that our standard monomials form a basis:  
first we show independence, then use the
Demazure  character formula to argue that our modules
have the same dimension as the number of standard
tableaux. Essential to the proof are the
vanishing theorems of 
Mathieu and Kumar for Bott-Samelson varieties
[Mt1,Mt2,Ku]. 
%In Section 4 we explain the connection with
%the work of Berenstein, Fomin, and Zelevinsky [BFZ]
%and use it to extend our results to unions of
%Bott-Samelson varieties.  

Reiner and Shimozono [RS2] give another
combinatorial interpretation of our tableaux.
There are also intriguing connections
between our basis and that of Brian Taylor [T]
for a special case of our modules.

\head 1. Definitions and main results \endhead

\subhead 1.1. Tableaux \endsubhead

We will use integer sequences to index several types
of objects.  Hence we will call an integer sequence
a ``word'', a ``tableau'', etc., depending on what
it indexes in a given context.

A {\it word} is a sequence $\ii = (i_1,\ldots,i_l)$
with $i_j \in \{1,\ldots,n-1\}$, and $\ii$ is 
{\it reduced} if $s_{i_1} \cdots s_{i_l} = w \in W$ is 
a minimal-length decomposition of $w$ into simple 
reflections. 

A {\it tableau} is a sequence
$\tau = (r_1,\ldots,r_N)$ with $r_j \in \{1,2,\ldots,n\}$.
For $\tau = (r_1,\ldots,r_N)$,\ $\tau'=(r_1',\ldots,r_{N'}')$,
we define the {\it concatenation}
$$
\tau \st \tau' = (r_1,\ldots,r_N,r_1',\ldots,r_{N'}').
$$
Let $\emp$ denote the empty tableau and define
$\emp \st \tau = \tau \st \emp = \tau$ for any
tableau $\tau$.

A {\it column} of size $i$ is a
tableau $\ka = (r_1,\ldots,r_i)$ with
$1 \leq r_1 < \cdots < r_i \leq n$.
The symmetric group $W$ acts on columns as follows:
for a permutation $w$ on $n$ letters
and a column $\ka = (r_1,\ldots,r_i)$,
the column $w\cdot \ka$ is the increasing rearrangement
of $(w(r_1),\ldots,w(r_i))$.
The {\it fundamental weight columns}
are the initial sequences:
$$
\om_i = (1,2,\ldots,i).
$$
The {\it Bruhat order} 
on columns is defined by elementwise comparison:  
$\ka = (r_1,\ldots,r_i) \leq \ka'=(r'_1,\ldots,r'_i)$
if and only if $r_1 \leq r'_1,\ldots,r_i\leq r'_i$.

For a word $\ii = (i_1,\ldots,i_l)$, and 
a sequence $\mm = (m_1,\ldots,m_l)$
with $m_j \in \ZZ^+$, we define a 
{\it tableau  of shape} $(\ii,\mm)$ to be
a concatenation of $m_1$ columns of size $i_1$,
\ $m_2$ columns of size $i_2$, etc:
$$
\tau = \ka_{11} \st \ka_{12} \st \cdots \st \ka_{1m_1}
\st \ka_{21} \st \cdots \st \ka_{2m_2} 
\st \cdots \st \ka_{l1} \st \cdots \st \ka_{lm_l},
$$
where $\ka_{km}$ is a column of size $i_k$ for each $k,m$.
(If $m_k=0$, there are {\it no} columns in the corresponding
position of $\tau$.)

\smallskip \noindent {\bf Remarks.}
(a) This terminology is suggested by the classical
notion of a column-strict Young tableau with $m_1$
columns of size $i_1$ followed by 
\  $m_2$ columns of size $i_2$,
etc., transcribed in terms of its column reading word. 
For example, take $\ii = (3,2,3)$,\  $\mm = (0,2,1)$,
which corresponds to the Young diagram at left below.
Note that $m_1 = 0$ means there are zero columns
of size $i_1 = 3$  on the left end of our diagram.
$$
\lam = 
\matrix 
\sq & \sq & \sq \\
\sq & \sq & \sq \\
    &     & \sq 
\endmatrix
\qquad \qquad
\tau = 
\matrix
3 & 1 & 1 \\
4 & 3 & 2 \\
  &   & 3
\endmatrix
$$
The filling at right is transcribed in our notation as
$\tau = 34 \st 13 \st 123$.
One can define a generalized Young diagram corresponding to
any reduced $(\ii,\mm)$.  (See [Mg2,Mg3,RS2].)
\newline
(b) In the model of Littelmann, bases are parametrized
by piecewise-linear paths in the weight lattice $\ZZ^n
\!\mod \ZZ (1,\ldots,1)$ of $G$.  
Our tableaux encode such paths if we consider
a column $(r_1,\ldots,r_i)$ as denoting a weight
$\pi = e_{r_1}+\cdots+e_{r_i}$ of the $i$th fundamental 
representation of $G$, so that a tableau is a sequence
of weights $\pi_1 \st \pi_2 \st \cdots$.  The associated
path goes in linear steps from from $0$ 
to $\pi_1$ to $\pi_1 + \pi_2$ etc.

\subhead 1.2 Liftable-standard tableaux \endsubhead
Let us once and for all {\bf arbitrarily fix  a  word}
$$
\ii = (i_1,\ldots,i_l),
$$
reduced or non-reduced.
From now on we will assume the presence of this chosen 
ambient word.
For $k \in \ZZ^+$, we will frequently use the notation
$$
[k] = \{1,2,\ldots,k\} ,
$$
as well as $[k,l]= \{k,k+1,\ldots,l\}$.

A {\it subword} of $\ii$ is a subsequence
$\ii' = (i_{j_1},i_{j_2},\ldots,i_{j_r})$
for some indices $1 \leq j_1 < \cdots < j_r\leq l$.
We say the set $J = \{j_1,\ldots,j_r\}\subseteq [l]$ 
is the {\it subword index} of $\ii'$, and 
we write $\ii'=\ii(J)$.  Note that we consider 
subwords $\ii(J_1)$, $\ii(J_2)$ to be different
whenever $J_1 \neq J_2$, so that there is a total
of $2^l$ distinct subwords.
Abusing notation, we will frequently identify
an indexing set $J\subseteq [l]$ with 
the corresponding subword $\ii(J)$ of our
fixed ambient word $\ii$, and we will
call $J$ itself a subword.
The intersection, union, and complement of two subwords
are defined in the obvious way in terms of their
indexing sets.
For $k\leq l$, the interval $[k] \subseteq [l]$ indexes
an initial subword of $\ii$.

Given any subword $J \subseteq [l]$, 
define $w(J)$, the
{\it permutation generated by $J$}, as
the partial product of $s_{i_1} s_{i_2} \cdots s_{i_l}$ 
containing only those factors which appear in 
$J=\{j_1<\cdots <j_r\}$:
$$
w(J) \, =\, \prod_{j \in J} s_{i_j}
\,=\, s_{i_{j_1}}\!\!\cdots s_{i_{j_r}}.
$$
Again, the subword $J$ is {\it reduced}
if the above is a minimal-length decomposition 
of $w(J)$ into $s_i$'s.
Also define the
{\it column generated by $J$ up to position $k$} to be
$$
w(\!J\!\cap\! [k]) \cdot \om_{i_k}.
$$

Now, consider a decreasing nest of subwords of $\ii$,
$$
[l] \supseteq J_{11} \supseteq \cdots \supseteq J_{1m_1}
\supseteq J_{21} \supseteq \cdots \supseteq J_{lm_l},
$$
We say it is a {\it reduced nest} if
$$
J_{km}\cap [k] \text{\quad is a reduced word for all } k,m.
$$
We say that a tableau $\tau$
of shape $(\ii,\mm)$ is
{\it generated}  by the reduced nest of subwords 
(or that the reduced nest 
is a {\it lifting} of the tableau)
if each column  $\ka_{km}$ of 
$\tau = \ka_{11}\!\st\! \cdots \!\st\! \ka_{lm_l}$ is generated 
by the subword $J_{km}$ up to the position $k$:
$$
\ka_{km} = w(\!J_{km}\!\cap\! [k]) \cdot \om_{i_k}.
$$
\medskip \noindent
{\bf Definition.}  A tableau $\tau$ of shape $(\ii,\mm)$
is {\it liftable-standard} (or just {\it standard})
if there exists a reduced nest of subwords of $\ii$ which
generates $\tau$.  The set of all standard tableaux
of shape $(\ii,\mm)$ is denoted $\T(\ii,\mm)$. 

A tableau $\tau$ is called {\it standard with 
respect to a subword} $J \subseteq [l]$ 
if all of the subwords  $J_{km}$
in the lifting are subwords of $J$:
$\ J \supseteq J_{11} \supseteq \cdots \supseteq J_{lm_l}$.
The set of such tableaux is denoted $\TT(J,\mm)$.
\medskip \noindent  
While useful to deduce general properties of standard tableaux,
this definition is quite difficult to work with in specific cases.
We will give a description of $\TT(\ii,\mm)$ 
allowing efficient computations in Section 1.4.

\medskip \noindent {\bf Examples.} (a)
Let $\ii = (1,2,1) = 121$,\ $\mm = (1,1,1)$ as in the 
introduction.
A typical subword index is $J = \{1,3\}$, associated
to the subword $\ii(J)=(i_1,i_3)=(1,1)$.
In order to emphasize that the position of the letters
is essential to distinguish subwords, we will write
$\miss$ in place of a letter of $\ii$ which is 
missing in $\ii(J)$.
That is, $\ii(J) = i_1\miss\, i_3 = 1\miss 1$.
For $J_1=\{1\}$, $J_2=\{3\}$, we have
$\ii(J_1)=1\miss\miss \neq \ii(J_2) = \miss\miss 1$,
and for the empty word we have
$\ii(\emptyset) = \miss\miss\miss$.

The nest of sets $J_{11}=\{1,2\}
\supseteq J_{21}=\{1,2\}\supseteq J_{31}=\{2\}$
indexes the nest of subwords $1 2 \miss
\supseteq 1 2 \miss \supseteq \miss 2 \miss$,
which generates the standard tableau
$\tau = s_1\om_1 \st s_1 s_2\om_2 \st s_2\om_1 
= 2\st 23 \st 1$. Another lifting for the same tableau is 
$1 2 1 \supseteq 1 2 1 \supseteq \miss\miss\miss$.
\smallskip \noindent
(b) Consider a $GL(n)$ Demazure 
module $V_w(\lam)$ for a permutation $w \in W$ and a 
partition $\lam = (\lam_1\geq \lam_2 \geq \cdots \geq \lam_n)$.
(The character of $V_w(\lam)$ is called a {\it key polynomial}.)\ \
Then $H(\Zii,\Lm)$ is isomorphic to the dual module $V_w^*(\lam)$ 
if $(\ii, \mm)$ are taken as follows.

Let $s_{i_1}\ldots s_{i_l} = w$ be a reduced decomposition.
Further suppose that if the last occurence of each letter 
$k =1,\ldots,n\!-\!1$ in $\ii$ is at position $j_k$, so that $i_{j_k} = k$,
then $j_1 < j_2 < \ldots <j_{n-1}$.
Now let $\lam' = (\lam'_1\geq \lam'_2 \geq \cdots \geq \lam'_n)$ be the
conjugate partition, and take $\mm = (m_1,\ldots,m_l)$ 
with $m_{j_k} = \lam'_k-\lam'_{k+1}$ and $m_j = 0$ otherwise.
That is, $m_{j_k}$ is the number of columns of size $k$ in 
the Young diagram of $\lam$.

It is easily seen that a classical Young 
tableau is semi-standard exactly if its column reading
word is liftable-standard with respect to the above
$(\ii, \mm)$ for $V_{w_0}(\lam)$, where $w_0$ is the longest
permutation.  The liftable-standard tableaux for the $(\ii,\mm)$
corresponding to a general $V_w(\lam)$ are exactly the
standard tableaux on the Schubert variety $X_w$ in classical
Standard Monomial Theory [LkSd1,LkSd2].

For example, the pair $\ii = (3,2,3)$, \ $\mm = (0,2,1)$ 
of the previous section give the Demazure
module $V_{s_3 s_2 s_3}(3,3,1)$.
The filling pictured is
a semi-standard Young tableau in the classical sense
and is standard on $X_{s_3 s_2 s_3}$.
Its column word 
$\tau =  34 \st 13\st 123$ has several liftings,
such as $3 2 \miss
\supseteq \miss 2 \miss \supseteq \miss\miss \miss$
and $3 2 3 \supseteq \miss 2 3 \supseteq 
\miss \miss 3$.  
\smallskip \noindent
(c) For any permutation $w \in W$ we have a kind of
generalized Young diagram called a Rothe diagram.
In [Mg3] we explain how to relate this to a pair $(\ii,\mm)$
so that $H^0(\ii,\mm)$ is the dual Schubert module of Kraskiewicz 
and Pragacz [KP], whose character is a Schubert polynomial.
In this case our standard tableaux are essentially identical
to the non-commutative Schubert polynomials of 
Lascoux and Schutzenberger [LcSb2].

\subhead 1.3. Standard basis \endsubhead

In the Introduction we defined the Bott-Samelson
variety $\Zii$, the embedding $\mu:\Zii \rightarrow \Gr(\ii)$,
and the line bundle $\Lm = \mu^* \OO(\mm)$.

Let 
$$
x=
\pmatrix
x_{11} & \cdots & x_{1i} \\
\vdots & \ddots & \vdots \\
x_{n1} & \cdots & x_{ni}
\endpmatrix
\in \Gr(i)
$$
be the homogeneous coordinates on the Grassmannian,
so that $x$ represents the subspace spanned by the 
column vectors of the matrix.
Then any column $\ka = (r_1 < \cdots < r_i)$
is associated to a Plucker coordinate, the
minor on rows $r_1,\ldots,r_n$ of $x$:
$$
\Del_{\ka}(x) \ = \
\det_{i\times i}
\pmatrix
x_{r_11} & \cdots & x_{r_1i} \\
\vdots & \ddots & \vdots \\
x_{r_i1} & \cdots & x_{r_ii}
\endpmatrix
\ \in \ H^0(\Gr(i),\OO(1)).
$$
Furthermore, the set of all tableaux of shape
$(\ii,\mm)$ parametrize a 
spanning set of $H^0(\Gr(\ii),\OO(\mm))$
consisting of monomials in the Plucker coordinates.
That is, for 
$\tau = \ka_{11} \st \cdots \st \ka_{lm_l}$, let
$$
\Del_{\tau} \, = \,
\prod_{j = 1}^l \prod_{m = 1}^{m_j} \Del_{\ka_{jm}}(x^{(j)})
\ \in\ H^0(\Gr(\ii),\OO(\mm)),
$$
where $x^{(j)}$ denotes the homogeneous coordinates
on the $j$th factor of $\Gr(\ii)$.
We let $\Del_{\emp} = 1$.
Denote the restriction of the section $\Del_{\tau}$ to
$\Zii \subseteq \Gr(\ii)$ by the same symbol $\Del_{\tau}$.
Under this restriction the Plucker monomials still form a spanning
set of $H^0(\Zii,\Lm)$ by the following ``Borel-Weil-Bott'' theorem:

\proclaim {Proposition (Mathieu {\rm [Mt1,Mt2]},
Kumar {\rm [Ku])}}
\newline 
(i) The map 
$$
\mu^* : H^0(\Gr(\ii),\OO(\mm)) \rightarrow
H^0(\Zii,\Lm)
$$ 
is a surjective homomorphism of 
$B$-modules. \newline
(ii) $H^i(\Zii,\Lm)=0$ for all $i > 0$.
\endproclaim

If a tableau $\tau$ is standard, we call $\Del_{\tau}$
a {\it standard monomial}.  

\proclaim {Theorem 1}  The standard monomials
of shape $(\ii,\mm)$ form a basis of the space of sections of 
$\Lm$ over $\Zii$:
$$
H^0(\Zii,\Lm) = 
\bigoplus_{\tau \in \T(\ii,\mm)}\!\!\! \FF\, \Del_{\tau}.
$$
\endproclaim \noindent
The proof will be given in Section 3.

Writing a diagonal matrix as $\diag(x_1,\ldots,x_n) \in T$,
we obtain the coordinate ring
$\FF[T]=\FF[x_1^{\pm 1},\ldots,x_n^{\pm 1}]$
(modulo the relation $x_1\!\cdots x_n\! = \! 1$ 
in case $G=SL(n)$\,).
By the (dual) character of a $B$-module $M$,
we mean 
$$
\Char M = \tr(\diag(x_1,\ldots,x_n)|M^*)\,
\in\, \FF[T].
$$
(We take duals in order to get polynomials
as characters.)
Now, given any tableau $\tau = (r_1,\ldots,r_N)$, 
we define its {\it weight monomial} 
$$
x^{\tau} = x_{r_1}\! \cdots\, x_{r_N} \in \FF[T].
$$
Then $\Char \FF \Del_{\tau} = x^{\tau}$, and we obtain:

\proclaim {Corollary}
$$
\Char H^0(\Zii,\Lm) = \sum_{\tau \in \T(\ii,\mm)} x^{\tau}.
$$
\endproclaim

\subhead 1.4. Demazure operations on tableaux \endsubhead

Define Demazure's isobaric divided difference
operator
$$
\Lam_i : \FF[T] 
\rightarrow \FF[T],
$$
$$
\Lam_i f = {x_i f - x_{i+1} s_i f \over x_i - x_{i+1}}.
$$

\smallskip \noindent {\bf Example.} Let
$f(x_1,x_2,x_3) = x_1^2 x_2^2 x_3$, so that
$$
\matrix
\Lam_2 f(x_1,x_2,x_3) &=&
{ x_2(x_1^2 x_2^2 x_3) - x_3(x_1^2 x_3^2 x_2) 
\over x_2 - x_3 } \\
&=& x_1^2 x_2 x_3 (x_2+x_3).
\endmatrix
$$

\medskip
%For any {\it reduced} decompostion of a permutation
%$w = s_{i_1}\ldots s_{i_l}$, define
%$$
%\Lam_{w} \eqdef \Lam_{i_1} \cdots \Lam_{i_l},
%$$
%which is known to be independent of the reduced
%decomposition chosen.

\noindent
Let $\om_i = x_1 x_2 \cdots x_i \in\FF[T]$, 
the $i$th fundamental weight of $GL(n)$.

\proclaim {Proposition (Demazure's formula {\rm [D2,Mt1,Ku]})}
$$
\Char H^0(\Zii,\Lm)  = 
\Lam_{i_1} ( \om_{i_1}^{m_1} \Lam_{i_2} 
( \om_{i_2}^{m_2} \cdots
\Lam_{i_l} (\om_{i_l}^{m_l}) \cdots )).
$$
\endproclaim

\medskip
Now we define analogs of the Demazure operations acting
on tableaux instead of on characters.  This will allow
us to ``lift'' the Demazure formula 
from characters to tableaux,  
thus reconciling the two character 
formulas above.  It also gives an efficient algorithm 
for generating the standard monomial basis.

We will need the {\it root operators} on tableaux 
first defined by Lascoux and Schutzenberger [LcSb1], 
and later generalized by Littelmann [Lt1, Lt2, Lt3].
For $i \in \{1,\ldots,n-1\}$,
the {\it lowering operator} $f_i$ takes a tableau 
$\tau=(r_1,r_2,\ldots)$ either 
to a formal null symbol $\nul$,
or to a new tableau 
$\tau' = (r'_1,r'_2,\ldots)$ 
by changing a single entry $r_j=i$ to $r'_j=i+1$
and leaving the other entries alone ($r'_j = r_j$),
according to the following rule.  

First, we ignore all the entries of 
$\tau$ except those equal to $i$ or $i+1$;  
if an $i$ is followed by an $i+1$ 
(not counting any ignored entries in between),
then henceforth we ignore that pair of entries; 
we look again for an $i$ followed 
(up to ignored entries) by an $i\!+\!1$,
and henceforth ignore this pair; 
and iterate until we ignore everything but 
a subsequence of the form 
$i\!+\!1,\, i\!+\!1,\ldots, i\!+\!1,\, i,\, i, \ldots, i$.
If there are {\it no} $i$ entries in this subsequence, 
then $f_i(\tau) = \nul$, the null symbol.  
If there {\it are} some $i$ entries,
then the {\it leftmost} is changed to $i\!+\!1$.

This is identical to 
Littelmann's minimum-point definition
[Lt1,Lt2] if we think of tableaux as paths
in the weight lattice.

\smallskip \noindent {\bf Example.}
We apply $f_2$ to the tableau
$$
\matrix
\tau & = & 1 & 2 & 2 & 1 & 3 & 2 & 1 & 4 & 2 & 2 & 3 & 3 \\
   &   & . & 2 & 2 & . & 3 & 2 & . & . & 2 & 2 & 3 & 3 \\
   &   & . & 2 & . & . & . & 2 & . & . & 2 & . & . & 3 \\
   &   & . & {\bold 2} & . & . & . & {\bold 2} & . & . & . & . & . & .
\\
\ \ f_2(\tau) & = &
 1 & {\bold 3} & 2 & 1 & 3 & {\bold 2} & 1 & 4 & 2 & 2 & 3 & 3  \\
(\!f_2\!)^2(\tau) & = &
 1 & {\bold 3} & 2 & 1 & 3 & {\bold 3} & 1 & 4 & 2 & 2 & 3 & 3  \\
(\!f_2\!)^3(\tau) & = & \nul
\endmatrix
$$
\medskip
\noindent
We also have the {\it raising operator} 
defined by $e_i(\tau) = (f_i)^{-1}(\tau)$ if this 
exists, $e_i(\tau) = \nul$ otherwise.
One can describe the action of $e_i$ identically
to that of $f_i$ except that $e_i$ changes 
the rightmost non-ignored $i\!+\!1$ into $i$.

Now define the plactic Demazure operator $\Lam_i$
taking a tableau $\tau$ to a set of tableaux:
$$
\Lam_i(\tau) = \{ \tau, f_i(\tau), f_i^2(\tau), \ldots \}
- \{\nul\}.
$$
To apply $\Lam_i$ to a set of tableaux $\T$, apply 
it to each element and take the union:
$$
\Lam_i(\T) = \bigcup_{\tau \in \T} \Lam_i(\tau).
$$

We will also need a tableau analog of multiplying
by a monomial in the $x_i$'s. For a column $\ka$, define
$\ka^{\st m} = \ka \st \cdots \st \ka$ ($m$ factors).
Then the multiplication by the monomial 
$\om_i^{m}$ in the character formula
will correspond to concatenating with the tableau
$\om_i^{\st m} = (1,2,\ldots,i,\ldots,1,2,\ldots,i)$.

Now we can build up the set of standard tableaux 
using the above operations. 

\proclaim {Theorem 2} The set of liftable-standard tableaux
is generated by the Demazure and concatenation operations:
$$
\T(\ii,\mm) = \Lam_{i_1}\left(\, \om_{i_1}^{\st m_1} \st 
\Lam_{i_2} \left(\, \om_{i_2}^{\st m_2} \st \cdots 
 \Lam_{i_l} \left(\om_{i_l}^{\st m_l}\right)\cdots  \right) \right).
$$
\endproclaim
\noindent The proof is given in Section 2.

\smallskip \noindent {\bf Example.} 
Let $\ii = 121$,\ $\mm=(1,1,1)$, so that 
$\TT(\ii,\mm) = \Lam_1( 1\! \st\! 
\Lam_2( 12\! \st\! \Lam_1( 1)))$.
To {\it generate} the standard tableaux,
we start with the empty tableau $\emptyset$,
and proceed from the right end of the above Demazure formula:
$$
\{ \emptyset \} {\overset 1\st  \to \rightarrow}
\{ 1 \} {\overset \Lam_1 \to \rightarrow}
\{ 1, 2 \} {\overset 12\st  \to \rightarrow}
\{ 12 1, 12 2 \} {\overset \Lam_2 \to \rightarrow}
$$
$$
\{ 12 1, 13 1, 12 2, 13 2, 13 3 \} {\overset 1\st  \to \rightarrow}
\{ 1 12 1, 1 13 1, 1 12 2, 1 13 2, 1 13 3 \} {\overset \Lam_1 \to \rightarrow}
$$
$$
\{ 1 12 1, 2 12 1, 2 12 2 ,
1 13 1, 2 13 1, 2 23 1, 2 23 2, 
1 12 2,
1 13 2, 2 13 2,
 1 13 3, 2 13 3, 2 23 3 \} 
$$ 
The last set is $\T(\ii,\mm)$.

To {\it test} whether a given tableau is standard,
say $\tau = 2123$, we invert the above operations:
at the $k$th step we raise the tableau as far as
possible using $f_{i_k}^{-1}$, then strip off
the initial word $\om_{i_k}$, then go on to the
next step.  That is,
$$
\tau = 2123 
{\overset f_1^{-1} \to\longrightarrow}
1123
{\overset (1\st)^{-1} \to\longrightarrow}
123
{\overset (12\st)^{-1} \to\longrightarrow}
3
$$
This algorithm will terminate in the empty tableau $\emp$
exactly if the the original $\tau$ is standard.
But in this case we end with a tableau $\tau' = 3$
for which we can invert neither $f_{i_3} = f_1$
nor $(\om_{i_3}\st) = (1\st)$, so the original
$\tau$ is {\it not} standard.  This process is
closely related to the keys of Lascoux and Schutzenberger
[LcSb1,LcSb2].

\medskip \noindent
The root operators $f_i$ and $e_i$ 
also define a crystal graph structure on $\T(\ii,\mm)$,
which suggests that our standard basis will deform to
a crystal basis inside the quantum function ring
of $B$.

\head  2. Proof of Theorem 2 \endhead

For a subword $J \subseteq [l]$, the set of 
{\it constructible tableaux} is 
$$
\matrix
\CC(J,\mm) &=& \Lam_{i_1}^{\del_1} \left(\, \om_{i_1}^{\st m_1} 
\st \Lam_{i_2}^{\del_2} \left(\, \om_{i_2}^{\st m_2} 
\st \cdots  \Lam_{i_l}^{\del_l}\left( 
\om_{i_l}^{\st m_l}\right)\cdots  \right) \right)\\
&=& \left\{ \ f_{i_1}^{a_1} ( \om_{i_1}^{\st m_1}\! \st\! f_{i_2}^{a_2} ( 
\om_{i_2}^{m_2}\! \st\! \cdots f_{i_l}^{a_l}(\om_{i_l}^{m_l}) \cdots )) 
\ \left| \
\matrix 
a_1,\ldots, a_l \geq 0,\\
 a_j = 0 \,\text{ for } j \not\in J
\endmatrix
\right. 
\ \right\}
\endmatrix
$$
where $\del_k = 1$ if $k \in J$, 
\ $\del_k = 0$ if $k \not \in J$.
We defined $\TT(J,\mm)$ in \S1.2.
(By convention, for $\mm = (0,\ldots,0)$ we set 
$\CC(J,\mm) = \TT(J,\mm) = \{\emptyset\}$, containing
only the empty tableau.)

We will show that 
$$
\TT(J,\mm) = \CC(J,\mm)
$$
for all subwords $J \subseteq [l]$.
We proceed by a series of elementary lemmas establishing 
identical recursions for the two sides.

\subhead 2.1 Recursion for $\TT(J,\mm)$ \endsubhead

\proclaim{Definition-Lemma 1} 
If $J \subseteq [l]$ is any subword 
and $J',J'' \subseteq J$
are maximal-length reduced subwords of $J$,
then $w(J') = w(J'')$.  We denote
$\wmax(J) \eqdef w(J')$ for any maximal
reduced $J' \subseteq J$.
\endproclaim 
\noindent
{\it Proof.}
By the subword definition of Bruhat order,
the lemma is equivalent to saying that 
$S(J) = \{w(J') \mid J'\subseteq J\}$
is an interval in the Bruhat order:
$S(J) = [e,\wmax]$
for some $\wmax \in W$.
By induction, we suppose this is true 
for $J$ and show it holds for 
$J_0 = \{j_0\} \cup J$, 
where $j_0 < j$ for all $j\in J$. 
Let $\wJ = \wmax(J)$. 
If $s_{j_0}\wJ > \wJ$,
then $S(J_0) = [e,s_{j_0} \wJ]$.  
If $s_{j_0} \wJ< \wJ$,
then $[e, s_{j_0} \wJ] \subseteq [e,\wJ]$,
and $S(J_0) = S(J)$, by the Zigzag Lemma [Hu \S 5.9].
\medskip 
\noindent
Note that if $J$ is reduced then $\wmax(J) = w(J)$.
\smallskip
Given $J' \subset J \subseteq [l]$, 
we say  $J'$ is {\it less than
$J$ with respect to column $k$ }
if the maximal $k$th column generated by 
$J'$ is smaller in Bruhat order than 
the maximal $k$th column generated by $J$:
$$
J' \lcol J \qquad \Leftrightarrow \qquad
\matrix  J' \subset J \quad \text{ and } \\
\wmax(J' \cap [k]) \cdot \om_{i_k} \ <\
\wmax(J \cap [k]) \cdot \om_{i_k}.
\endmatrix
$$
Now, let $\ep(k) = (0,\ldots,1,\ldots,0)$, a sequence of
length $l$ with a 1 in the $k$th place.
Then for $\mm=(0,\ldots,0,m_k,\ldots, m_l)$,
we have $\mm-\ep(k) = (0,\ldots,0,m_k-1,\ldots,m_l)$.
\proclaim{Lemma 2}
For $\mm = (0,\ldots,0,m_k,\ldots,m_l)$ with $m_k>0$,
$J \subseteq [l]$, and \newline
$\kmax = \wmax(J\cap [k]) \cdot \om_{i_k}$,
we have
$$
\TT(J,\mm) = \kmax \st \TT(J,\mm-\ep(k)) \ \,
\sqcup
\bigcup_{J' \lcol J }
\! \TT(J',\mm).
$$
\endproclaim
\noindent
{\it Proof.}
(a) First, it is evident that 
$\TT(J',\mm) \subseteq \TT(J,\mm)$ for any $J' \subseteq J$.
\newline
(b)  Also, $\kmax \st \TT(J,\mm-\ep(k)) \subseteq \TT(J,\mm)$
as follows.
If $\tau' = \ka_1 \st \ka_2 \st \cdots
\in \TT(J,\mm-\ep(k))$,
by definition there exists a lifting 
$J_1 \supseteq J_2 \supseteq \cdots$ 
with $J_1 \subseteq J$ and $J_1 \cap [k]$ reduced.

Now let $\tJ \subseteq [k]$ be a maximal reduced
subword of $J\cap [k]$. By the Definition-Lemma,
$$
w(\tJ) \geq w(J_1 \cap [k])\geq w(J_2 \cap [k])\geq \cdots,
$$ 
so we may take reduced words $\tJ_k \subseteq \tJ$ with 
$w(\tJ_j) = w(J_j \cap [k])$ and 
$\tJ_1 \supseteq \tJ_2 \supseteq \cdots$.
Set $J'_j = J_j\! \cap [k\!+\!1,l]$.
Then
$$
\tJ \!\cup\! [k\!+\!1,l]\, \supseteq \,\tJ_1\! \cup\! J'_1
\,\supseteq\, \tJ_2 \!\cup\! J'_2 \supseteq \cdots
$$
is a lifting of 
$\kmax \st \tau' = \kmax \st \ka_1 \st \ka_2 \cdots$.
\newline
(c)  Finally, suppose 
$\tau = \ka_0 \st \ka_1\st \ka_2 \cdots \in \TT(J,\mm)$, 
with lifting $J_0 \supseteq J_1 \supseteq J_2 \supseteq \cdots$.
Then we must have exactly one of the following. 
Either $\ka_0 = \kmax$ and 
$J_1 \supseteq J_2 \supseteq \cdots$ is
a lifting of $\ka_1 \st \ka_2 \st \cdots$, so that 
$\tau \in \ka \st \TT(J,\mm-\ep(k))$.
Or $\ka_0 \neq \kmax$, meaning
$$
\wmax(J_0\cap [k]) \cdot \om_{i_k}
\neq \wmax(J \cap [k]) \cdot \om_{i_k},
$$
and hence 
$$
\wmax(J_0\cap [k]) \cdot \om_{i_k}
< \wmax(J \cap [k]) \cdot \om_{i_k}.
$$
Therefore $J_0 \lcol J$ and $\tau \in \TT(J_0,\mm)$.
\proclaim{Lemma 3}
Let $\mm = (0,\ldots,0,m_k,\ldots,m_l)$.
If $J,J' \subseteq [l]$ are subwords with
$J \cap [k\!+\!1,l] = J' \cap [k\!+\!1,l]$
and $\wmax(J\cap [k]) = \wmax(J'\cap [k])$,
then $\TT(J,\mm) = \TT(J',\mm)$. 
\endproclaim
\noindent {\it Proof.}
If $\ka_1 \st \ka_2\cdots \in \TT(J,\mm)$ 
has a lifting $J_1 \supseteq J_2 \cdots$ with $J_1 \subseteq J$,
then 
$$
\wmax(J\cap [k]) \geq w(J_1 \cap [k]) \geq w(J_2 \cap [k]) \geq \cdots,
$$
so by the subword definition of Bruhat order
we can find  reduced words with 
$$
J'\cap [k]\supseteq \tJ_1 \supseteq \tJ_2 \supseteq \cdots 
$$ 
with $w(\tJ_j) = w(J_j\cap [k])$.
Setting $J'_j = \tJ_j \cup (J_j \cap [k\!+\!1,l])$ for all $j$
gives a lifting $J'_1 \supseteq  J'_2 \supseteq \cdots$ for
$\ka_1 \st \ka_2 \st \cdots$ which shows it to lie in $\TT(J',\mm)$.

Reversing the roles of $J$ and $J'$ we obtain the reverse inclusion,
which completes the proof.

\subhead 2.2. Head-string property \endsubhead

\medskip\noindent
{\bf Definition.}
For $i \in \{1,\ldots,n-1\}$, 
the {\it $i$-string through $\tau$} is defined as
$$
S_i(\tau) \eqdef \{\cdots, e_i^2 \tau, e_i \tau, \tau, 
f_i \tau, f_i^2 \tau, \cdots\}-\{\nul\}.
$$
We say a set of tableaux $\TT$ has the 
{\it head-string property} 
if for any $i$ and any $\tau \in \TT$,
we have either : \newline
(i) $S_i(\tau) \subseteq \TT$ (the entire $i$-string of
$\tau$ lies in $\TT$); or \newline
(ii) $S_i(\tau) \cap \TT = \{\tau\}$ and
$e_i\tau = \nul$ (only the head of the
string lies in $\TT$). \medskip \noindent
A key step in our proof of Theorems 1 and 2 will be to 
show that the set of standard tableaux has this property. 

We will use the following properties
special to groups of type $A$:  for any $i$ and any column $\ka$,
$$
e_i\,\ka = \nul \, \text{ or } \, f_i\,\ka = \nul
\qquad \text{and} \qquad
e_i^2\,\ka = f_i^2\,\ka = \nul.
$$

\proclaim{Lemma 4}
For $\ka$ a column, $\tau'$ a tableau, and $a > 0$,
we have
$$
f_i^a (\ka \st  \tau') =
\left\{
\matrix
(f_i \kappa) \st (f_i^{a-1}\tau')
& \text{ if } f_i\ka \neq \nul,\ e_i\tau' = \nul \\
\kappa \st (f_i^a \tau') &
\text{ otherwise, }
\endmatrix
\right.
$$
$$
e_i^a (\ka \st  \tau') =
\left\{
\matrix
\kappa \st (e_i^{a}\tau') &
 \text{ if } 
(f_i\ka = \nul \text{ and } e_i\tau' \neq \nul)\\
&\text{ or } e_i^2 \tau' \neq \nul \\
(e_i \kappa) \st (e_i^{a-1}\tau') &
\text{ otherwise. }
\endmatrix
\right.
$$
Here we use the convention that $\tau \st \nul = \nul \st \tau = \nul$.
\endproclaim
\noindent
{\it Proof.} This follows from the well-known 
(and easily checked) formulas [Lt2]
$$
f_i(\tau \st \tau') =
\left\{
\matrix
(f_i \tau) \st \tau'
& \text{ if } \forsome n>0, \ f_i^n \tau \neq \nul,\ e_i^n \tau' = \nul\ \\
\tau \st (f_i\tau') &
\text{ otherwise, }
\endmatrix
\right.
$$
$$
e_i(\tau \st \tau') =
\left\{
\matrix
\tau \st (e_i\tau') &
\text{ if }\forsome n>0, \ f_i^n \tau = \nul,\ e_i^n \tau' \neq \nul\ \\
(e_i \tau) \st \tau' & 
\text{ otherwise,}
\endmatrix
\right.
$$
together with $f_i^2\kappa = \nul$.

\subhead 2.3. Recursion for $\CC(J,\mm)$ \endsubhead

\proclaim{Theorem $\bold{2}^+$\!} \newline
(i) As in Lemma 2, 
let $\mm = (0,\ldots,0,m_k,\ldots,m_l)$ with $m_k>0$,
$J \subseteq [l]$, and \newline
$\kmax = \wmax(J\cap [k]) \cdot \om_{i_k}$.
Then we have
$$
\CC(J,\mm) = \kmax \st \CC(J,\mm-\ep(k)) \ \,
\sqcup
\bigcup_{J' \lcol J }
\! \CC(J',\mm).
$$
(ii) $\CC(J,\mm) = \TT(J,\mm)$. 
\smallskip
\noindent
(iii)  $\CC(J,\mm)$\ \ has the head-string property.
\endproclaim
\noindent
{\it Proof.} 
By induction on $|J|$ (the order of $J$) and $|\mm| = m_1+\cdots+m_l$
(the number of columns in a tableau).  The initial cases
$J = \emptyset$ or $\mm=(0,\ldots,0)$ are trivial.
Now assume (i)-(iii) for all
$\CC(J',\mm')$ with $|J'|<|J|$ or $|\mm'|<|\mm|$.
We will prove (i)-(iii) for $\CC(J,\mm)$.
\medskip \noindent
(i) If $J \cap [k] = \emptyset$, the righthand side of the
equation (i) reduces to $\om_{i_k} \st \CC(J,\mm-\ep(k))$,
and the claim is clear.

Otherwise, let $j_1$ be the smallest element of $J\cap [k]$, and write
$$\tJ = J - \{j_1\},\qquad i = i_{j_1},\qquad 
\tkmax = \wmax(\tJ \cap [k])\! \cdot\! \om_{i_k}\ \leq \ \kmax.
$$
Then $\CC(J,\mm) = \Lam_{i} \CC(\tJ,\mm)$, so that 
every tableau in $\CC(J,\mm)$ may be written as
$\ka \st \tau = f_i^a(\tka\st \ttau)$ for 
$\tka \st\ttau \in \CC(\tJ,\mm)$.
In fact, we have  
$\tau \in \CC(J,\mm-\ep(k))$, which follows easily by  
induction and Lemma 4.

By induction, we have
$$
\CC(\tJ,\mm) = \tkmax \st \CC(\tJ,\mm-\ep(k)) \ \,
\sqcup
\bigcup_{\tJ' \lcol \tJ }
\! \CC(\tJ',\mm).
$$
(a)  First, it is evident that $\CC(J,\mm) \supseteq \CC(J',\mm)$
whenever $J \supseteq J'$. 
\newline
(b)  We show $\CC(J,\mm) \supseteq \kmax \st \CC(J,\mm-\ep(k))$
as follows. \newline
For $\tau \in \CC(J,\mm-\ep(k))$, we can write
$\tau = f_i^a \ttau$ for $\ttau \in \CC(\tJ,\mm-\ep(k))$.
In fact, by raising $\ttau$ with $e_i$ and increasing $a$, 
we may assume 
that $e_i \ttau = \nul$.  
(We know $e_i^a \ttau \in \CC(\tJ,\mm-\ep(k)) \cup \{\nul\}$
by the head-string property (iii).)

In  case  $\kmax = \tkmax$, we have $f_i \tkmax = \nul$, and
$$
\kmax\! \st\! \tau \,=\, \kmax\! \st\! (f_i^a \ttau) 
\,=\, f_i^a(\tkmax\! \st\! \ttau)
\ \in\ \Lam_i \CC(\tJ,\mm) \,=\, \CC(J,\mm).
$$

In case $\kmax > \tkmax$, we have $\kmax = f_i \tkmax$, and 
$$
\kmax \!\st\! \tau = \kmax \!\st\! (f_i^a \ttau) 
= f_i^{a+1}(\tkmax \!\st\! \ttau) 
\ \in\ \Lam_i \CC(\tJ,\mm) \,=\, \CC(J,\mm).
$$
This completes the $\supseteq$ direction of formula (i).
\newline
(c)  Now we show the $\subseteq$ direction of formula (i).
%\newline
We suppose $\ka \st \tau \in \CC(J,\mm)$  with $\ka < \kmax$,
and proceed to show $\ka \st \tau \in \CC(J',\mm)$ as in 
the Theorem.  Let us write 
$\ka \st \tau  = f_i^a(\tka \st \ttau)$ with
$\ttau \in \CC(\tJ,\mm-\ep(k))$.

In case $\kmax = \tkmax$ we have $\tka<\tkmax$, so by (i)
applied to $\CC(\tJ,\mm)$,\  
$\tka \st \ttau \in \CC(\tJ',\mm)$ for some $\tJ' \lcol \tJ$.
We may assume $\tka = \tkmax' 
\eqdef \wmax(\tJ' \cap [k])\cdot \om_{i_k}$,
since otherwise we would have $\tka \st \ttau$ in some smaller
$\CC(\tJ'',\mm)$ by induction.
If $f_i \tkmax' < \kmax$, then by definition 
$\{j_1\} \cup\! \tJ' \lcol J$,
so $\ka\st \tau \in \CC(\{j_1\} \cup \tJ',\mm)$ 
gives the desired result.
If $f_i \tkmax' = \kmax$, then (since $\ka < \kmax$) we must have
$f_i^a(\tkmax' \st \ttau) = \tkmax' \st (f_i^a \ttau)$, and 
therefore  $e_i \ttau \neq \nul$ by Lemma 4.  
This means $f_i^a \ttau \in \CC(\tJ',\mm-\ep(k))$
by head-string, and so 
$$
\ka \st \tau = \tkmax' \st (f_i^a \ttau) \in \CC(\tJ',\mm)
$$ 
by (i) applied to $\CC(\tJ',\mm)$.  Since $\tJ' \lcol J$,
we have the desired result.

In case $\kmax > \tkmax$, we have $\tJ \lcol J$.
If $\tkmax > \tka$, 
then $\tkmax \st \ttau \in \CC(\tJ',\mm)$
for $\tJ' \lcol \tJ$, so that 
$\{j_1 \} \cup \tJ' \lcol J$,
and clearly $\ka \st \tau \in \CC(\{j_1 \} \cup \tJ',\mm)$, 
as desired.
If $\tkmax = \tka$, we must have 
(since $\ka < \kmax$ and $f_i\tkmax = \kmax$) that
$f_i^a(\tkmax\st\ttau) = \tkmax\st (f_i^a \ttau)$, which means
$e_i \ttau \neq \nul$ by Lemma 4.  Thus, by head-string, we have
$f_i^a \ttau \in \CC(\tJ,\mm-\ep(k))$, so that, as desired,
$$
\ka \st \tau \,=\, \tkmax \st (f_i^a\ttau) \ \in \ 
\tkmax \st \CC(\tJ,\mm-\ep(k))
\,\subseteq\, \CC(\tJ,\mm).
$$

This completes the proof of (i).
\smallskip \noindent
(ii)  Follows immediately from the identical recursions
for $\TT(J,\mm)$ (Lemma 2) and for $\CC(J,\mm)$ (part (i)),
by induction on $|J|$ and $|\mm|$.
The intitial cases $J=\emptyset$ and $\mm=(0,\ldots,0)$ are
trivial. 
\smallskip \noindent
(iii)  To show the head-string property for $\CC(J,\mm)$,
we need to prove:  for all $i_0$,
$$
\tau' \in \CC(J,\mm)\ \text{ and } \
e_{i_0} \tau' \neq \nul \quad \Rightarrow \quad 
e_{i_0} \tau' \in \CC(J,\mm)\ \text{ and } 
f_i \tau' \in \CC(J,\mm) \cup \{\nul\}.
$$
Take $\tau'=\ka \st \tau$.  
If $\ka < \kmax$, then by (i) we have
$\tau' \in \CC(J',\mm)$ with $J'\lcol J$, and
the head-string property follows by induction.

Thus we may assume $\tau' = \kmax \st \tau$ 
with $\tau\in\CC(J,\mm-\ep(k))$.

In case $e_{i_0} \kmax = \nul$, we have by hypothesis
$e_{i_0} (\kmax \st \tau) \neq \nul$, so we must have
$e_{i_0} (\kmax \st \tau) = \kmax \st (e_{i_0} \tau)$.
Hence $e_{i_0} \tau \neq \nul$, and by head-string
$e_{i_0} \tau \in \CC(J,\mm-\ep(k))$, and
$$
\kmax \st e_{i_0} \tau\ \in\ \kmax \st \CC(J,\mm-\ep(k))
\ \subseteq\  \CC(J,\mm)
$$
by (i), as desired.  Also by head-string $f_{i_0} \tau \in \CC(J,\mm-\ep(k))
\,\cup \{\nul\}$, and 
$$
\kmax \st f_{i_0} \tau\ \in\ \kmax \st \CC(J,\mm-\ep(k))\, \cup \{\nul\}
\ \subseteq\  \CC(J,\mm)\, \cup \{\nul\}
$$
as desired.

In case $f_{i_0}\kmax = \nul$, 
if $e_{i_0}(\kmax \st \tau) = \kmax \st (e_{i_0} \tau) \neq \nul$,
we may argue as in the previous case.
Thus suppose $e_{i_0}(\kmax \st \tau) = 
(e_{i_0} \kmax) \st \tau \neq \nul$, 
so that $e_{i_0} \kmax \neq \nul$.
Now let $\wmax = \wmax(J\cap [k])$, so that 
$\kmax = \wmax\cdot \om_{i_k}$.
Since $f_{i_0} \kmax = \nul$,\ $e_{i_0} \kmax \neq \nul$,
we have $s_{i_0} \kmax = e_{i_0}\kmax < \kmax$,
and so $s_{i_0} \wmax < \wmax$.  Therefore we may 
find  a reduced word for $\wmax$ with first letter equal to $i_0$:
$$
\tii = (i'_1,\ldots,i'_t)  \quad \text{ with } \quad
w(\tii) = \wmax\ \ \text{and}\ \ i'_1 = i_0.
$$
Let 
$$
\ii' = (i'_1,\ldots,i'_t,i'_{t+1},i'_{t+2},\ldots)
\quad \text{ with } \quad
(i'_{t+1},i'_{t+2},\ldots) = \ii(J\cap [k\!+\!1,l]).
$$
Using (ii) and Lemma 3, we have
$$
\CC(J,\mm)=\TT(J,\mm)=\TT(\ii',\mm)=\CC(\ii',\mm)
$$
(Note that $|\,\ii'|\leq |J|$, so (ii) holds for $\CC(\ii',\mm)$.)
\ \ But $\CC(\ii',\mm) = \Lam_{i_0}(\cdots)$, so
for any $\tau' \in \CC(\ii',\mm)$, we have
$e_{i_0} \tau', f_{i_0} \tau'\in \CC(\ii',\mm) \cup \{\nul\}
= \CC(J,\mm) \cup \{\nul\}$ as desired.

The proof of (iii) is finished, the induction proceeds,
and the Theorem is proved.

\head 3. Proof of Theorem 1 \endhead

\subhead 3.1 Subvarieties \endsubhead

Given a subword index $J \subseteq [l]$ we may 
consider $\ii(J)$ as a word in its own right,
corresponding to a Bott-Samelson variety
$Z_{\ii(J)}$ which embeds naturally into
$\Zii$ via
$$
Z_{\ii(J)} \cong Z_J \eqdef Q_1 \times \cdots \times Q_l \ / \ B^l
\ \subseteq \ \Zii = P_{i_1} \times \cdots \times P_{i_l} 
\ /\ B^l,
$$
where 
$$
Q_j = \left\{
\matrix
P_{i_j} & \text{ if } j \in J \\
B & \text{ if } j \not\in J
\endmatrix \right. .
$$

Let us index Schubert varieties $X_{\ka}$
in a Grassmannian $\Gr(i)$
by columns $\ka  = (r_1,\cdots,r_i)$.
That is, let $\FF^{\ka}$ be the subspace of $\FF^n$
spanned by the coordinate vectors $e_r$ for $r \in \kappa$,
and define $X_{\ka} = \overline{B \cdot \FF^{\ka}} 
\subseteq \Gr(i)$.
(Under the isomorphism $\Gr(i) \cong G/P$ for a suitable
maximal parabolic $P$, we can 
write this as $X_{\kappa} 
\cong \overline{B\cdot w P} \subseteq G/P$,
where $\ka=w \cdot \om_i$.)
We have $X_{\ka} \subseteq X_{\ka'}$ if and
only if $\ka \leq \ka'$ in Bruhat order.

\proclaim{Lemma} For any $J \subseteq [l]$,
the partial multiplication map 
$$
\matrix 
\mu_k : & Z_J & \rightarrow &
\Gr(i_k) \\
&(p_1,p_2,\ldots,p_l) & \mapsto & p_1\!\!\cdots\! p_k\FF^{i_k}
\endmatrix
$$
has image equal to the Schubert variety of 
the column generated by $J$ up to position $k$:
$$
\roman{Im}(\mu) = X_{\kappa}, \qquad
\kappa = \wmax(J \cap [k]) \cdot \om_{i_k}.
$$
\endproclaim
\noindent
{\it Proof.} This follows from the formula:
$$
P_i X_{\ka} = 
\left\{
\matrix
X_{s_i \ka} & \text{ if }  s_i \ka >\ka \\
X_{\ka} & \text{ otherwise.}
\endmatrix \right. 
$$

We denote the restriction of the line bundle $\Lm$ 
from $\Zii$ to $Z_J$ by the same symbol $\Lm$.
In order to prove Theorem 1, we will show the more
general fact that
$\TT(J,\mm)$ indexes a basis of $H^0(Z_J,\Lm)$.

\subhead 3.2 Linear independence \endsubhead

For any subwords $J_1, J_2, \cdots \subseteq [l]$,
we may consider the union of the corresponding
Bott-Samelson
varieties  embedded in $\Gr(\ii)$:
$$
Z_{J_1} \cup Z_{J_2} \cup \cdots \subseteq \Gr(\ii).
$$
The restriction of $\OO(\mm)$ again defines
a line bundle $\Lm$ on the union, and 
for any tableau $\tau$ the Plucker monomial
$\Del_{\tau}$ restricts to an element of
$H^0(Z_{J_1} \cup Z_{J_2} \cup \cdots, \Lm)$.
\smallskip \noindent
{\bf Definition.}  A tableau $\tau$ of shape
$(\ii,\mm)$ is 
{\it standard on a union} 
$Z = Z_{J_1} \cup Z_{J_2} \cup \cdots\ $  
if it is standard on at least one of the components 
$Z_{J_1}, Z_{J_2},\ldots$.
That is, the set of standard tableaux on $Z$ is 
$$
\TT(Z_{J_1}\cup Z_{J_2} \cup \cdots,\mm)
\eqdef \TT(J_1,\mm) \cup \TT(J_2,\mm) \cup \cdots.
$$ 

\proclaim {Proposition}
For any subwords $J_1, J_2,\ldots \subseteq [l]$,
the standard monomials of shape $(\ii,\mm)$ on 
the union $Z = Z_{J_1} \cup Z_{J_2} \cup \cdots$ are
linearly independent.
\endproclaim

\noindent {\bf Remark.}  A statement of this generality
holds only for independence:  the standard monomials
on a union $Z$ do {\it not} in general span $H^0(Z,\Lm)$.

For example, let $\ii = (1,2,1)$,\ $\mm = (0,0,1)$,\
$J_1 =\{1\}$,\ $J_2 =\{3\}$.
Then $\TT(J_1,\mm) = \TT(J_2,\mm) = \{1,2\}$,\ 
but $\dim H^0(Z_{J_1}\! \cup\! Z_{J_2},\Lm) = 3$.
In fact, in this case the restriction
map $H^0(\Gr(\ii),\OO(\mm)) 
\rightarrow H^0(Z_{J_1} \!\!\cup\!\! Z_{J_2},\Lm)$ is not
surjective.  This is possible because $\OO(\mm)$
is non-ample.

\medskip
\noindent
{\it Proof of Proposition.} Let 
$\tau^{(1)},\ldots,\tau^{(t)}$ be standard tableaux
in $\TT(Z,\mm)$.  Consider a linear relation among
the standard monomials $\Del_{\tau^{(r)}}$ on the
variety Z
$$
(*) \qquad \qquad
 a_1 \Del_{\tau^{(1)}}(z) + \cdots + 
a_t \Del_{\tau^{(t)}}(z) = 0
\quad\ \forall\, z\! \in\! Z, \qquad \qquad \qquad
$$
where $a_r \in \FF$.  We will show
$$
a_r =0 \text{ \quad for \ \ } r = 1,\ldots,t
$$
by induction on $t$ (the length of the linear relation)
and on $|\mm| = m_1+\cdots+m_l$ (the number of columns
in a tableau). \newline
(a) Let us suppose $\mm = (0,\ldots,0,m_k,\ldots,m_l)$
with $m_k >0$, and write $\tau^{(r)} = 
\ka_{k1}^{(r)}\st \cdots \st\ka_{lm_l}^{(r)}$.
Let $I_{k1}^{(r)} \supseteq \cdots 
\supseteq I_{lm_l}^{(r)}$
be a lifting of $\tau^{(r)}$. 
By definition, each $I_{k1}^{(r)}$ is contained
in one of the subwords $J_1,J_2,\ldots$ defining $Z$,
so the Bott-Samelson variety of the subword $I_{k1}^{(r)}$
is contained in $Z$:
$$
Z_{I_{k1}^{(r)}} \subseteq Z \text{\quad for all } r.
$$
(b) Now let $\ka$ denote one of the Bruhat-minimal
elements among the first columns of 
$\tau^{(1)}, \ldots,\tau^{(t)}$:
$$
\ka \in \min\{\, \ka_{k1}^{(1)},\ldots,\ka_{k1}^{(t)}\ \}
$$
Order the terms of relation $\,(*)\,$ so that,
for some $1 \leq t_0 \leq t$, we have
$$
\ka = \ka_{k1}^{(r)} \ \ \ \text{for \ } r \leq t_0,
\qquad \qquad 
\ka \not \geq \ka_{k1}^{(r)} \ \ \ \text{for \ } r > t_0.
$$
(c) We show that $a_r=0$ for $r \leq t_0$.
Let
$$
Y = \bigcup_{r \leq t_0}\!\! Z_{I_{k1}^{(r)}} \ \subset \ Z.
$$
Let us restrict the relation \ $(*)$\ from $Z$ to $Y$.
By the Lemma of \S3.1, we have
$ \mu_k(Y) = X_{\ka} \subseteq \Gr(i_k)$.  Furthermore,
the first factor $\Del_{\ka_{11}^{(r)}}(z)$ 
of $\Del_{\tau^{(r)}}$ is just the Plucker coordinate of
$\ka_{k1}^{(r)}$ on $\Gr(i_k)$.  Since 
$\ka \not\geq \ka_{k1}^{(r)}$ for all $r>t_0$,
we have 
$$
\Del_{\tau^{(r)}}(y) = 0 \quad \forall\, y \in Y,
\text{\qquad for all }\ r >t_0.
$$
so that $\ (*)\ $ becomes
$$
\Del_{\ka}(y) \ 
\left(\, a_1 \Del_{\ttau^{(1)}}(y)+ \cdots 
+ a_{t_0}\Del_{\ttau^{(t_0)}}(y) \,\right) =0 
\qquad \forall y\! \in\! Y,
$$
where $\tau^{(r)} = \ka_{k1}^{(r)} \st \ttau^{(r)}$
for some $\ttau^{(r)} \in \TT(I_{k1}^{(r)},\mm-\ep(k))$.
However by the same Lemma,
$\Del_{\ka}$ is not identically zero on any of
the components $Z_{I_{k1}^{(r)}}$ of $Y$.  
Hence $\Del_{\ka}$ is not
a zero-divisor in the coordinate ring of $Y$,
and we may factor it from the equation
to get a linear relation among
standard monomials $\ttau^{(r)}$ of shape 
$(\ii,\mm-\ep(k))$ on $Y$.  That is,
$\ttau^{(r)} \in \TT(Y,\mm-\ep(k))$, and 
$$
a_1 \Del_{\ttau^{(1)}}(y)+ \cdots 
+ a_{t_0}\Del_{\ttau^{(t_0)}}(y) =0 
\qquad \forall\, y\! \in\! Y.
$$
By induction on $|\mm|$, this relation must be
identically zero: $a_r = 0$ for $r \leq t_0$.
\newline
(d)  Since $t_0 \geq 1$,
we have shown that $a_r=0$ for at least a single $r$.
Therefore $\ (*)\ $ reduces to a relation with
fewer than $t$ terms, which must have $a_r = 0$ for all $r$
by induction on $t$.  The proof of the Proposition is finished.

\subhead 3.3  Dimension counting \endsubhead

Recall that for $\tau = (r_1,r_2,\ldots)$ we define
$x^{\tau} = x_{r_1} x_{r_2} \cdots \in \FF[T]$.
For any set of tableaux $\TT$, let 
$$
\Char \TT \eqdef \sum_{\tau \in \TT} x^{\tau},
$$

\proclaim{Proposition}  
For any subword $J\subseteq [l]$ and $\mm=(m_1,\ldots,m_l)$,
$m_j \geq 0$, we have
$$
\Char \TT(J,\mm)= 
\Lam_{i_1}^{\delta_1} ( \om_{i_1}^{m_1} \Lam_{i_2}^{\delta_2} 
( \om_{i_2}^{m_2} \cdots
\Lam_{i_l}^{\delta_l} (\om_{i_l}^{m_l}) \cdots )),
$$
where $\delta_j = 1$ if $j\in J$ and
$\delta_j =0$ if $j\not\in J$.
\endproclaim
\medskip \noindent
{\it Proof.}  By Theorem $2^+$ of \S 2.3, we know that 
$$
\TT(J,\mm)= 
\Lam_{i_1}^{\delta_1} ( \om_{i_1}^{\st m_1}\!\st\! \Lam_{i_2}^{\delta_2} 
( \om_{i_2}^{\st m_2}\!\st\! \cdots
\Lam_{i_l}^{\delta_l} (\om_{i_l}^{\st m_l}) \cdots )),
$$
so we need to show that each operation
$\Lam_i$ and $(\om_i\st)$ on tableaux has the
corresponding effect on characters.

For $j \in [l]$, let 
$$\mm \cap [j,l] = (0,\ldots,0,m_{j},\ldots,m_l),
\qquad \TT_{j}=\TT(J\cap [j,l],\, \mm\cap [j,l]).
$$
Now, 
$$
\om_{i_k}^{\st m_k}\st \TT_{k+1} =
\TT(J\,\cap\, [k\!+\!1,l],\mm\,\cap\, [k,l]),
$$
so by \S 2.3 this set has the head-string property.
That is, we may partition it into 
$i_k$-strings
$$
\om_{i_k}^{\st m_k}\st \TT_{k+1}
= S^{(1)} \sqcup S^{(2)} \sqcup \cdots
$$
so that each $S^{(r)}$ is either a complete $i_k$-string
or only the head of an $i_k$-string.  It is easily
verified that
$$
\Char( \Lam_{i_k} S^{(r)}) = \Lam_{i_k}(\Char S^{(r)}),
$$
so we have
$$
\matrix
\Char \TT_k 
& = & \Char\, \Lam_{i_k}(\om_{i_k}^{\st m_k}\st \TT_{k+1})\,  \\
& = & \Char\, \Lam_{i_k}(S^{(1)} \sqcup S^{(2)} \sqcup \cdots ) \\
& = &  \Lam_{i_k}\,\Char(S^{(1)} \sqcup S^{(2)} \sqcup \cdots ) \\
& = & \Lam_{i_k}\,\Char(\om_{i_k}^{\st m_k}\st \TT_{k+1})\,  \\
& = & \Lam_{i_k}\,(\,\om_{i_k}^{m_k}\,\, \Char\TT_{k+1} ) 
\endmatrix
$$
Thus we may build up $\TT(\ii,\mm)$ and $\Char \TT(\ii,\mm)$
in parallel steps, and the Proposition follows.
\medskip
\noindent
{\bf Proof of Theorem 1.}
The Demazure character formula of \S 1.4 applies to the subvarieties
$Z_J$ to  give
$$
\Char H^0(Z_J,\Lm)= 
\Lam_{i_1}^{\delta_1} ( \om_{i_1}^{m_1} \Lam_{i_2}^{\delta_2} 
( \om_{i_2}^{m_2} \cdots
\Lam_{i_l}^{\delta_l} (\om_{i_l}^{m_l}) \cdots )).
$$
Hence by the above Proposition we have
(after specializing the characters to $x_1=\cdots=x_n=1$):
$$
\#\{\Del_{\tau}\mid \tau \in \TT(J,\mm)\}
= \dim H^0(Z_J,\Lm).
$$
But by the Proposition of \S 3.2,
we know that the standard monomials
$\{\Del_{\tau}\mid \tau \in \TT(J,\mm)\}$
form a linearly independent subset of $H^0(Z_J,\Lm)$.
Therefore they form a basis.

\enddocument